\documentclass[11pt]{article}
\usepackage{amsmath, amssymb, amscd, epsfig}

\newtheorem{thm}{Theorem}
\newtheorem{prop}[thm]{Proposition}
\newtheorem{lem}[thm]{Lemma}

\newtheorem{cor}[thm]{Corollary}
\newtheorem{rem}[thm]{Remark}
\newtheorem{df}[thm]{Definition}
\newtheorem{ex}[thm]{Example}

\renewcommand{\epsilon}{\varepsilon}
\renewcommand{\phi}{\varphi}

\newcommand{\BB}{\mathbb}
\newcommand{\g}{\mathfrak}
\newcommand{\pf}{\noindent {\it Proof. }}
\newcommand{\qed}{\nopagebreak $\qquad$ $\square$ \vskip5pt}
\newcommand{\separate}{\vskip5pt}
\newcommand{\VF}{\operatorname{VF}}

\newcommand{\supp}{\operatorname{supp}}
\newcommand{\im}{\operatorname{Im}}
\newcommand{\re}{\operatorname{Re}}

\begin{document}

\title{\bf Integrals of equivariant forms over non-compact symplectic
manifolds}
\author{Matvei Libine}
\maketitle

\begin{abstract}
This article is a result of the AIM workshop on Moment Maps and
Surjectivity in Various Geometries (August 9 -- 13, 2004)
organized by T.~Holm, E.~Lerman and S.~Tolman.
At that workshop I was introduced to the work of
T.~Hausel and N.~Proudfoot on hyperk\"ahler quotients \cite{HP}.
One interesting feature of their article is that they consider integrals
of equivariant forms over non-compact symplectic manifolds which do not
converge, so they formally {\em define} these integrals as sums over
the zeroes of vector fields, as in the Berline-Vergne localization formula.
In this article we introduce a geometric-analytic regularization technique
which makes such integrals converge and utilizes the symplectic structure of
the manifold. We also prove that the Berline-Vergne localization formula
holds for these integrals as well.
The key step here is to redefine the collection of integrals
$\int_M \alpha(X)$, $X \in \g g$, as a distribution on the Lie algebra $\g g$.
We expect our regularization technique to generalize to non-compact group
actions extending the results of \cite{L1,L2}.
\end{abstract}


%\tableofcontents

\begin{section}
{Introduction}
\end{section}

This article is a result of the AIM workshop on Moment Maps and
Surjectivity in Various Geometries (August 9 -- 13, 2004)
organized by T.~Holm, E.~Lerman and S.~Tolman.
At that workshop I was introduced to the work of
T.~Hausel and N.~Proudfoot on hyperk\"ahler quotients \cite{HP}.
One interesting feature of their article is that they consider integrals
of equivariant forms over non-compact symplectic manifolds which do not
converge, so they formally {\em define} these integrals as sums over
the zeroes of vector fields, as in the Berline-Vergne localization formula.

While the definition is perfectly valid, it does not feel satisfactory.
The Berline-Vergne localization formula relates a global object (integral
of a cohomology class) with a local object (certain quotients defined
at zeroes of a vector field). From this point of view, the localization
formula is very similar to the Lefschetz fixed point formula.
The Lefschetz fixed point formula fails for non-compact manifolds in general
because the fixed points may ``run away to infinity''.
By analogy one expects the Berline-Vergne formula to fail on non-compact
manifolds for the same reason.

In this article we introduce a geometric-analytic regularization technique
for non-compact symplectic manifolds which makes such integrals converge
in the subanalytic setting.
Then we prove under an additional assumption that ensures the zeroes do not
run away to infinity (see Definition \ref{no-zeroes}) that the Berline-Vergne
localization formula holds for these integrals as well.
The key step here is to redefine the collection of integrals
$\int_M \alpha(X)$, $X \in \g g$, as a distribution on the Lie algebra $\g g$.
Such approach has been used before, particularly in the important special
case when the manifold $M$ is a coadjoint orbit of a real semisimple Lie group
(see, for example, \cite{DHV}, \cite{DV}, \cite{KV}, \cite{L1,L2},
\cite{Par1, Par2} and references therein).
We extend these ideas to a much wider class of subanalytic symplectic manifolds.
The subanalytic setting allows us to use the properties of
o-minimal structures (for elementary introductions see \cite{DM}, \cite{D})
to estimate the growth of integrals at infinity.
The subanalytic setting is more general than real algebraic setting
and, in fact, all results of this paper hold in even more general setting
of polynomially bounded o-minimal structures.
To make this article more accessible we discuss the most relevant properties of
o-minimal structures in the Appendix.

The main results of this paper are
Proposition \ref{subanalytic-hamiltonian-prop},
the regularization technique described in Subsection \ref{regularization}
and Theorems \ref{loc1}, \ref{main}.
Proposition \ref{subanalytic-hamiltonian-prop} has appeared in a very special
case as Proposition 8 in \cite{DV}.
Theorem \ref{loc1} is closely related to another extension of the
Berline-Vergne localization formula to non-compact setting given in \cite{L2}.
It is also related to the extension of the Duistermaat-Heckman formula
to non-compact manifolds due to E.~Prato and S.~Wu \cite{PrWu}.
Theorem \ref{main} essentially says that the Berline-Vergne localization
formula still holds after our regularization procedure.
As an example, we consider coadjoint orbits of real semisimple Lie groups.

It is interesting to note that a similar problem was studied by physicists
G.~Moore, N.~Nekrasov and S.~Shatashvili in \cite{MNS}.
They regularize volumes of non-compact k\"ahler and hyperk\"ahler quotients
by introducing a parameter $\epsilon$, discarding the poles in $\epsilon$
and extracting the finite part as $\epsilon$ goes to zero.

We expect our regularization technique to generalize to non-compact group
actions. An analogue of Berline-Vergne integral localization formula for
non-compact group actions was proved in \cite{L1,L2}. However, in \cite{L1,L2}
integrals are taken over homology cycles sitting inside a cotangent bundle
$T^*M$ which comes with the standard Hamiltonian structure,
the integrand has to satisfy many conditions as well.
In the subanalytic setting these results can be extended further to more
general manifolds and differential forms.

Finally, I would like to thank the reviewers for their comments, suggestions
and pointing out the reference \cite{MNS}.

\separate

\begin{section}
{Notations and definitions}
\end{section}

\subsection
{Equivariant cohomology}

Let $G$ be a compact real algebraic Lie group;
we denote by $\g g$ its Lie algebra and
by $\g g^*$ the dual of the Lie algebra.
The group $G$ acts on $\g g$ by the adjoint action and on $\g g^*$
by the coadjoint action.
Let $M$ be a smooth manifold which is possibly non-compact,
and let $G$ act on $M$.
%We use $\Omega^*(M)$ to denote the algebra of complex-valued
%differential forms on $M$.
The action of $G$ on $M$ lifts to an action on
$\Omega^*(M)$ -- the algebra of complex-valued
differential forms on $M$.
To each element $X \in \g g$ we associate a vector field
$\VF_X$ on $M$ defined by
$$
(\VF_X \cdot f) (m) =
\frac d{d\epsilon} f \bigl( \exp (\epsilon X)m \bigr) \Bigr|_{\epsilon=0},
\qquad f \in {\cal C}^{\infty}(M), \quad m \in M.
$$
We denote by $\iota_X$ the contraction by the vector field $\VF_X$, and
by $M^X$ the set of zeroes of $\VF_X$:
$$
M^X = \{ m \in M;\: \VF_X(m) =0 \}.
$$

For introductions to equivariant forms and equivariant cohomology see
\cite{BGV}, \cite{GS}, \cite{L3}.
Recall that a (smooth) $G$-equivariant form on $M$ is a ${\cal C}^{\infty}$-map
$\alpha: \g g \to \Omega^*(M)$ which is $G$-equivariant, i.e.
$$
\alpha(g \cdot X) = g \cdot (\alpha (X))
\qquad \text{for all $g \in G$ and $X \in \g g$},
$$
and possibly non-homogeneous.
We denote by $\Omega_G^{\infty}(M)$ the space of $G$-equivariant forms on $M$
and by $\Omega_G^*(M)$ the subalgebra of $\Omega_G^{\infty}(M)$ consisting of
equivariant forms depending on $X \in \g g$ polynomially.
We define a {\em twisted deRham differential} by
\begin{equation}  \label{d_g}
(d_{\g g} \alpha) (X) = d(\alpha(X)) + \iota_X (\alpha(X)),
\qquad X \in \g g,
\end{equation}
where $d$ denotes the ordinary deRham differential.
The map $d_{\g g}$ preserves $G$-equivariant forms, and $(d_{\g g})^2=0$ on
$\Omega_G^{\infty}(M)$.
An equivariant form $\alpha$ such that $d_{\g g}\alpha=0$ is called
{\em equivariantly closed}.
The equivariant degree of a polynomial equivariant form in $\Omega_G^*(M)$
is its differential form degree plus twice its degree as a polynomial on
$\g g$, then the twisted deRham differential $d_{\g g}$ increases the
equivariant degree of $\alpha$ by one.
On the other hand, for an element $\alpha \in \Omega_G^{\infty}(M)$,
the differential $d_{\g g}$ changes the parity of the differential
form degree of $\alpha$. We denote by $H_G^*(M)$ and by $H_G^{\infty}(M)$
the cohomologies of complexes $(\Omega_G^*(M), d_{\g g})$ and
$(\Omega_G^{\infty}(M), d_{\g g})$ respectively;
$H_G^*(M)$ is $\BB Z$-graded and called the (ordinary) $G$-equivariant
cohomology of $M$, $H_G^{\infty}(M)$ is $\BB Z_2$-graded and called the
$G$-equivariant cohomology of $M$ with smooth coefficients.

%If $\alpha$ is an equivariant differential form, $\alpha_{[k]}$ denotes
%its homogeneous differential form component of degree $k$.

When $\gamma \in \Omega^*(M)$ is a differential form and $k=0,1,2,\dots$,
we denote by $\gamma_{[k]}$ the homogeneous differential form component
of degree $k$. If $N \subset M$ is an oriented submanifold we define
$$
\int_N \gamma \quad =_{def} \quad
\int_N \bigl(\gamma_{[\dim N]} \bigr) \bigr|_N.
$$

\subsection
{Equivariant cohomology with distributional coefficients}

In this article we will be primarily interested in $G$-equivariant forms
and cohomology on $M$ with distributional coefficients.
Let $\Omega_c^{top}(\g g)$ the space of smooth compactly supported
complex-valued differential forms on $\g g$ of top degree;
it will play the role of the space of test functions.
The adjoint action of $G$ on $\g g$ lifts to an action on
$\Omega_c^{top}(\g g)$.
We equip both $\Omega_c^{top}(\g g)$ and $\Omega^*(M)$ with
${\cal C}^{\infty}$ topologies.
By an equivariant form with distributional (or ${\cal C}^{-\infty}$)
coefficients we mean a continuous $\BB C$-linear $G$-equivariant map
$$
\alpha: \quad \Omega_c^{top}(\g g) \ni \quad \phi \mapsto
\langle \alpha, \phi \rangle_{\g g} \quad \in \Omega^*(M),
$$
we denote the space of those by $\Omega_G^{-\infty}(M)$.
We treat elements $\alpha \in \Omega_G^{-\infty}(M)$ as $\Omega^*(M)$-valued
$G$-equivariant distributions on $\g g$.

To define the twisted deRham differential on $\Omega_G^{-\infty}(M)$, we
pick a real vector space basis $\{E^1,\dots,E^{\dim \g g}\}$ of $\g g$ and
let $\{E_1,\dots,E_{\dim \g g}\}$ be the associated dual basis of $\g g^*$.
We regard each $E_k$, $1 \le k \le \dim \g g$, as a linear function on $\g g$;
in particular the product $E_k \phi$ makes sense for
$\phi \in \Omega_c^{top}(\g g)$.
For each $\alpha \in \Omega_G^{-\infty}(M)$, we set
\begin{equation}  \label{d_g2}
\langle d_{\g g} \alpha, \phi \rangle_{\g g}=
d \langle \alpha, \phi \rangle_{\g g} +
\sum_{k=1}^{\dim \g g} \iota_{E^k} \langle \alpha, E_k \phi \rangle_{\g g},
\qquad \phi \in \Omega_c^{top}(\g g).
\end{equation}
We have a natural inclusion $\Omega_G^{\infty}(M) \subset \Omega_G^{-\infty}(M)$,
and the formulas (\ref{d_g}) and (\ref{d_g2}) agree on $\Omega_G^{\infty}(M)$.
The differential $d_{\g g}$ defined by (\ref{d_g2}) is independent of a
particular choice of the basis $\{E^1,\dots,E^{\dim \g g}\}$ of $\g g$, and
we still have $(d_{\g g})^2=0$ on $\Omega_G^{-\infty}(M)$.
The cohomology of $(\Omega_G^{-\infty}(M), d_{\g g})$ is denoted by
$H_G^{-\infty}(M)$, it is $\BB Z_2$-graded and called the
$G$-equivariant cohomology of $M$ with distributional (or ${\cal C}^{-\infty}$)
coefficients.

Let ${\cal C}_c^{\infty}(\g g^*)$ be the space of test functions, i.e.
the space of smooth compactly supported complex-valued functions on $\g g^*$
endowed with ${\cal C}^{\infty}$ topology.
We also consider the space of continuous $G$-equivariant $\BB C$-linear maps
from ${\cal C}_c^{\infty}(\g g^*)$ to $\Omega^*(M)$, denoted by
${\cal M}^{-\infty}(\g g^*, M)^G$.
%${\cal M}^{-\infty}(\g g^*, M)^G$ -- the subspace of such $G$-equivariant maps.
Similarly, we treat elements $\beta \in {\cal M}^{-\infty}(\g g^*,M)^G$
as $\Omega^*(M)$-valued $G$-equivariant distributions on $\g g^*$
and denote the value of $\beta$ at $\psi \in {\cal C}_c^{\infty}(\g g^*)$
by $\langle \beta, \psi \rangle_{\g g^*} \in \Omega^*(M)$.
Next we define the differential on ${\cal M}^{-\infty}(\g g^*, M)^G$.
For each $\beta \in {\cal M}^{-\infty}(\g g^*, M)^G$, we set
\begin{equation}  \label{d_g3}
\langle \widehat{d_{\g g}} \beta , \psi \rangle_{\g g^*} =
d \langle \beta, \psi \rangle_{\g g^*} - i \sum_{k=1}^{\dim \g g}
\iota_{E^k} \langle \beta, \partial_{E_k} \psi \rangle_{\g g^*},
\qquad \psi \in {\cal C}_c^{\infty}(\g g^*),
\end{equation}
where $\partial_{E_k} \psi$ denotes the partial derivative of $\psi$
relative to the basis $\{E_1,\dots,E_{\dim \g g}\}$ of $\g g^*$.
As before, the differential $\widehat{d_{\g g}}$ is independent of a
particular choice of the basis $\{E^1,\dots,E^{\dim \g g}\}$ of $\g g$, and
$(\widehat{d_{\g g}})^2=0$ on ${\cal M}^{-\infty}(\g g^*, M)^G$.

The complexes $(\Omega_G^{-\infty}(M), d_{\g g})$ and
$({\cal M}^{-\infty}(\g g^*, M)^G, \widehat{d_{\g g}})$ are related by a
Fourier transform.
We denote by
$$
\Omega_{temp}^{-\infty}(\g g,M)^G \subset \Omega_G^{-\infty}(M)
\qquad \text{and} \qquad
{\cal M}_{temp}^{-\infty}(\g g^*, M)^G \subset {\cal M}^{-\infty}(\g g^*,M)^G
$$
the subspaces of tempered $G$-equivariant distributions on
$\g g$ and $\g g^*$ respectively with values in $\Omega^*(M)$.
%Note that if $\phi \in \Omega_c^{top}(\g g)$, then
%$\int_{\g g} e^{i\langle \xi,X \rangle} \phi(X)$ is a rapidly decreasing
%function in $\xi \in \g g^*$.
The Fourier transform
${\cal F} : \Omega_{temp}^{-\infty}(\g g,M)^G \to
{\cal M}_{temp}^{-\infty}(\g g^*,M)^G$
is normalized so that
$$
\Bigl\langle {\cal F}(\alpha), \int_{\g g} e^{i\langle \xi,X \rangle} \phi(X)
\Bigr\rangle_{\g g^*} = \langle \alpha, \phi \rangle_{\g g},
\qquad \forall \alpha \in \Omega_{temp}^{-\infty}(\g g,M)^G ,\quad
\forall \phi \in \Omega_c^{top}(\g g).
$$
%which can be expressed as
%$$
%\int_{\g g^*} e^{i \langle \xi, X \rangle} {\cal F}(\alpha) (\xi) = \alpha(X)
%\qquad X \in \g g ,\: \xi \in \g g^*.
%$$
%(We do not need $G$-equivariance to define the Fourier transform,
%but we will not be interested in non-equivariant distributions anyway.)
Then the following diagram commutes:
$$
\begin{CD}
\Omega_{temp}^{-\infty}(\g g,M)^G @>{d_{\g g}}>>
\Omega_{temp}^{-\infty}(\g g,M)^G \\
@V{\cal F}VV      @VV{\cal F}V \\
{\cal M}_{temp}^{-\infty}(\g g^*,M)^G   @>{\widehat{d_{\g g}}}>>
{\cal M}_{temp}^{-\infty}(\g g^*,M)^G
\end{CD}
$$
i.e. the Fourier transform ${\cal F}$ becomes a chain map.

Following \cite{Par2}, we say that a distribution
$\beta \in {\cal M}^{-\infty}(\g g^*,M)^G$
has compact support in {\em $\g g^*$-mean} on $M$ if, for every test function
$\psi \in {\cal C}_c^{\infty}(\g g^*)$, the form
$\langle \beta, \psi \rangle_{\g g^*} \in \Omega^*(M)$ has compact support
in $M$.

\begin{df}
For every $\beta \in {\cal M}^{-\infty}(\g g^*,M)^G$ with
compact support in $\g g^*$-mean on $M$ we denote by $\int_M^{distrib} \beta$
the distribution on $\g g^*$ defined by
$$
{\cal C}_c^{\infty}(\g g^*) \ni \quad \psi \mapsto
\Bigl\langle \int_M^{distrib} \beta, \psi \Bigr\rangle_{\g g^*} =_{def}
\int_M \langle \beta, \psi \rangle_{\g g^*} \quad \in \BB C.
$$
\end{df}

\begin{lem}[Lemma 2.11 in \cite{Par2}] 
If $\beta \in {\cal M}^{-\infty}(\g g^*,M)^G$ has compact support
in $\g g^*$-mean on $M$, then so does $\widehat{d_{\g g}} \beta$, moreover
$$
\int_M^{distrib} \widehat{d_{\g g}} \beta =0.
$$
\end{lem}

%\pf
%Pick any $\psi \in {\cal C}_c^{\infty}(\g g^*)$, then (\ref{d_g3})
%implies that
%$$
%\supp \langle \widehat{d_{\g g}}\beta, \psi \rangle_{\g g^*}
%\quad \subset \quad
%\supp \langle \beta, \psi \rangle_{\g g^*} \cup \bigcup_{k=1}^{\dim \g g}
%\supp \langle \beta, \partial_{E_k} \psi \rangle_{\g g^*} \quad \subset M,
%$$
%hence $\widehat{d_{\g g}} \beta$ has compact support in $\g g^*$-mean too.
%We have
%$$
%\Bigl\langle \int_M^{distrib} \widehat{d_{\g g}} \beta,
%\psi \Bigr\rangle_{\g g^*} =
%\int_M \langle \widehat{d_{\g g}} \beta, \psi \rangle_{\g g^*} =
%\int_M d \langle \beta, \psi \rangle_{\g g^*} - i \sum_{k=1}^{\dim \g g}
%\int_M \iota_{E^k} \langle \beta, \partial_{E_k} \psi \rangle_{\g g^*}.
%$$
%The first term of the RHS is zero because it is an integral of an exact form,
%and each integral
%$\int_M \iota_{E^k} \langle \beta, \partial_{E_k} \psi \rangle_{\g g^*}$
%is also zero because the differential forms
%$\iota_{E^k} \langle \beta, \partial_{E_k} \psi \rangle_{\g g^*}$ have no
%components of maximal degree.
%\qed

\subsection
{Hamiltonian systems}

In this article we will be interested in the case when $M$ is symplectic
and the action of $G$ is Hamiltonian. That is, the manifold $M$ comes
equipped with a symplectic form $\omega$, the action of $G$ preserves $\omega$,
and there is a moment map $\mu: M \to \g g^*$ which is $G$-equivariant:
$$
\mu(g \cdot m) = g \cdot \mu (m)
\qquad \text{for all $g \in G$ and $m \in M$,}
$$
and such that
\begin{equation}  \label{hamiltonian}
d\mu = -\iota_X (\omega), \qquad \text{for all $X \in \g g$.}
\end{equation}
A symplectic manifold $(M, \omega)$ always has a preferred orientation
-- the one given by $\omega^{\frac 12 \dim M}$.

We can also regard $\mu$ as a linear map from $\g g$ to $\Omega^0(M)$ --
the space of functions on $M$.
Define $\tilde\omega \in \Omega_G^{\infty}(M)$ by
$$
\tilde\omega : \g g \to \Omega^0(M) \oplus \Omega^2(M) \subset \Omega^*(M),
\qquad \tilde \omega(X) = \mu(X) + \omega, \qquad X \in \g g.
$$
Then $\tilde \omega$ is equivariantly closed, which follows immediately
from (\ref{hamiltonian}).

Throughout this article we always assume that the moment map 
$\mu: M \to \g g^*$ is {\em proper},
i.e.
$$
\forall K \subset \g g^*, \text{ $K$ compact } \quad \Longrightarrow \quad
\mu^{-1} (K) \subset M \text { is compact}.
$$

\subsection
{P.-E.~Paradan's pushforward map
${\cal P}: H^*_G(M) \to {\cal C}^{-\infty}(\g g^*)$}

Let $(M,\omega,\mu)$ be a Hamiltonian system.
As always, we assume that the moment map $\mu$ is proper.
Pick an element $\alpha \in \Omega_G^{\infty}(M)$. Then, for each
$\phi \in \Omega_c^{top}(\g g)$, the integral
\begin{equation}  \label{my}
\int_{\g g} \alpha(X) \wedge e^{i\tilde\omega(X)} \wedge \phi(X)
= e^{i\omega} \wedge \int_{\g g} \alpha(X) \wedge e^{i\mu(X)} \wedge \phi(X),
\qquad X \in \g g,
\end{equation}
produces a differential form on $M$.
Thus $\alpha \wedge e^{i\tilde\omega}$ may be regarded as a $G$-equivariant
distribution on $\g g$ with values in $\Omega^*(M)$, i.e. an element of
$\Omega^{-\infty}(\g g,M)^G$.
We call this distribution $(\alpha \wedge e^{i\tilde\omega})^{distrib}$.

\begin{lem} [Lemma 2.12 in \cite{Par2}]
For every equivariant form $\alpha(X) \in \Omega_G^*(M)$ with polynomial
dependence on $X \in \g g$, the $\Omega^*(M)$-valued distribution
$$
(\alpha \wedge e^{i\tilde\omega})^{distrib}: \quad
\Omega_c^{top}(\g g) \ni \quad \phi(X) \mapsto
\int_{\g g} \alpha(X) \wedge e^{i\tilde\omega(X)} \wedge \phi(X)
\quad \in \Omega^*(M)
$$
is $G$-equivariant and tempered, hence belongs to
$\Omega_{temp}^{-\infty}(\g g,M)^G$.
Its Fourier transform ${\cal F}(\alpha \wedge e^{i\tilde\omega})^{distrib}$
has compact support in $\g g^*$-mean on $M$.
The distribution
$\int_M^{distrib} {\cal F}(\alpha \wedge e^{i\tilde\omega})^{distrib}$
is $G$-invariant.
If $\alpha$ is equivariantly exact, i.e. $\alpha = d_{\g g} \alpha'$
for some $\alpha' \in \Omega_G^*(M)$, then
$\int_M^{distrib} {\cal F}(\alpha \wedge e^{i\tilde\omega})=0$.
\end{lem}

P.-E.~Paradan's key observation is that, for every polynomial
$P(X)$ on $\g g$, we can associate a differential operator
$P(-i\partial_{\xi})$ (with constant coefficients) on
${\cal C}^{\infty}(\g g^*)$ so that
$$
P(X) \cdot {\cal F} \bigl( \psi(\xi) \bigr) =
{\cal F } \bigl( P(-i\partial_{\xi}) \psi(\xi) \bigr),
\qquad X\in \g g ,\quad \xi \in \g g^*,
$$
%$$
%P(X) \int_{\g g^*} e^{-i \langle \xi,X\rangle} \psi(\xi) \,d\xi =
%\int_{\g g^*} e^{-i \langle \xi,X\rangle}
%\bigl[P(-i\partial_{\xi}) \psi \bigr](\xi) \,d\xi,
%\qquad X\in \g g ,\: \xi \in \g g^*,
%$$
for all test functions $\psi(\xi) \in {\cal C}_c^{\infty}(\g g^*)$.
This association extends naturally to $\Omega^*(M)$-valued polynomials on
$\g g$ and hence to $\Omega_G^*(M)$; for an $\alpha(X) \in \Omega_G^*(M)$
we denote by $\alpha(-i\partial_{\xi})$ the corresponding differential
operator with values in $\Omega^*(M)$.
Then, for a test function $\psi \in {\cal C}_c^{\infty}(\g g^*)$
with $\hat\psi(X) = {\cal F}(\psi)$, we have:
\begin{multline}  \label{paradan}
\langle {\cal F}(\alpha \wedge e^{i\tilde\omega})^{distrib},
\psi \rangle_{\g g^*} \bigr|_m =
\int_{\g g} \alpha(X) \wedge e^{i\tilde\omega(X)} \wedge \hat\psi(X)
\biggr|_m  \\
= e^{i\tilde\omega} \bigl[\alpha(-i\partial_{\xi}) \psi \bigr](\mu(m)),
\qquad \forall m \in M.
\end{multline}
Hence the differential form
$\langle {\cal F}(\alpha \wedge e^{i\tilde\omega})^{distrib},
\psi \rangle_{\g g^*}$
is supported inside $\mu^{-1}(\supp \psi)$ which is compact because $\mu$
is proper.
This implies that ${\cal F}(\alpha \wedge e^{i\tilde\omega})^{distrib}$
has compact support in $\g g^*$-mean on $M$.

Using this lemma, P.-E.~Paradan defines a pushforward map
from $\Omega_G^*(M)$ into ${\cal C}^{-\infty}(\g g^*)$
-- distributions on $\g g^*$ by
$\alpha \mapsto \int_M^{distrib} {\cal F}(\alpha \wedge e^{i\tilde\omega})^{distrib}$
which descends to a map
$$
{\cal P}: H^*_G(M) \to {\cal C}^{-\infty}(\g g^*)
$$
(equation (14) in \cite{Par2}).

\begin{ex}
Let the equivariant form $\alpha=1$. Then
${\cal P}(1) = \int_M^{distrib} {\cal F}(e^{i\tilde\omega})^{distrib}$ is (up to a
constant coefficient) just the Duistermaat-Heckman measure on $\g g^*$.
(The Duistermaat-Heckman measure on $\g g^*$ is the pushforward of
the Liouville measure on $M$ to $\g g^*$ via the moment map.)
\end{ex}

In his paper \cite{Par2} P.-E.~Paradan studies the distribution
$\int_M^{distrib} {\cal F}(\alpha \wedge e^{i\tilde\omega})^{distrib}$ on $\g g^*$
and gives a localization formula for it at (the connected components of) the
critical points of $\|\mu\|_{\g g^*}^2$.

\subsection
{Integrals of equivariant forms}

In this article we will be primarily interested in the distribution
$(\alpha \wedge e^{i\tilde\omega})^{distrib}$ defined by (\ref{my}).
%Distributions of this kind have appeared in \cite{L1,L2}.
While
$\bigl\langle (\alpha \wedge e^{i\tilde\omega})^{distrib}, \phi \bigr\rangle_{\g g}
\in \Omega^*(M)$
need not have compact support, we will show that in the subanalytic setting
(or, even more generally, in the setting of polynomially bounded o-minimal
structures) this form is automatically integrable.
Hence we obtain a distribution which we denote by
$\int_M^{distrib} \alpha \wedge e^{i\tilde\omega}$ on $\g g$.
Comparing (\ref{my}) and (\ref{paradan}) we see that the distributions
$\int_M^{distrib} \alpha \wedge e^{i\tilde\omega}$ on $\g g$ and
${\cal P}(\alpha)$ on $\g g^*$ are related to each other
by the Fourier transform.
In particular, P.-E.~Paradan's results immediately apply
to $\int_M^{distrib} \alpha \wedge e^{i\tilde\omega}$ as well.

%Let $(M,\omega,\mu)$ be a Hamiltonian system.
%As always, we assume that the moment map $\mu$ is proper.
%Let $n = \frac 12 \dim M$.
We fix a positive definite inner product
$(\cdot,\cdot)_{\g g^*}$ on $\g g^*$ which is invariant
under the coadjoint action of $G$, and denote by $\|\cdot\|_{\g g^*}$
the corresponding norm. Let
$$
M_{\le R} = \{ m \in M;\: \|\mu(m)\|_{\g g^*} \le R\}, \qquad R>0.
$$
Since the moment map $\mu$ is proper, the sets $M_{\le R}$ are compact.
Note also that the inner product $(\cdot,\cdot)_{\g g^*}$ being $G$-invariant
implies that the sets $M_{\le R}$ are preserved by the $G$-action on $M$.

\begin{df}
We call an element $\alpha \in \Omega_G^{\infty}(M)$
{\em integrable in distributional sense} if the limit
$$
\lim_{R \to +\infty} 
\int_{M_{\le R}} \Bigl|
\int_{\g g} \alpha \wedge e^{i\tilde\omega} \wedge \phi \Bigr|
$$
exists for every test form $\phi \in \Omega^{top}_c(\g g)$ and the mapping
\begin{equation}   \label{int_M}
\Omega^{top}_c(\g g) \to \BB C,
\quad \phi \mapsto \lim_{R \to +\infty} 
\int_{M_{\le R}} \Bigl(
\int_{\g g} \alpha \wedge e^{i\tilde\omega} \wedge \phi \Bigr)
\end{equation}
is continuous.

For an integrable $\alpha$, we denote the distribution (\ref{int_M}) by
$\int_M^{distrib} \alpha \wedge e^{i\tilde\omega}$.
\end{df}

\begin{ex}  \label{ex1}
Suppose that the equivariant form $\alpha =1$ is integrable in distributional
sense.
Then the distribution $\int_M^{distrib} e^{i\tilde\omega}$ is
(up to a constant coefficient) the Fourier transform of the
Duistermaat-Heckman measure on $\g g^*$.
\end{ex}

We first show that in the subanalytic setting every subanalytic
$\alpha \in \Omega_G^*(M)$ is integrable in distributional sense,
then we study properties of the distribution 
$\int_M^{distrib} \alpha \wedge e^{i\tilde\omega}$.
%Subanalytic sets and functions are defined in  Definition \ref{subanalytic};
%for subanalytic differential forms see discussion at the end of
%Subsection \ref{structure-intro}.

\begin{df} \label{subanalytic-hamiltonian-def}
By a {\em subanalytic Hamiltonian system} we mean a Hamiltonian system
$(M,\omega,\mu)$ such that the manifold $M$, the group $G$, the symplectic
form $\omega$, the moment map $\mu$ and the action of the group $G$ on $M$
are subanalytic in the sense of Definition \ref{subanalytic} and
discussion at the end of Subsection \ref{structure-intro}.
\end{df}

\begin{prop}  \label{subanalytic-hamiltonian-prop}
Let $(M,\omega,\mu)$ be a subanalytic Hamiltonian system with proper
moment map $\mu$, and let $\alpha: \g g \to \Omega^*(M)$ be a subanalytic
form depending on $X \in \g g$ polynomially, then $\alpha$ is integrable
in distributional sense.
%i.e. the limit
%$$
%\lim_{R \to +\infty} 
%\int_{M_{\le R}} \Bigl|
%\int_{\g g} \alpha \wedge e^{i\tilde\omega} \wedge \phi \Bigr|,
%\qquad M_{\le R} = \{ m \in M;\: \|\mu(m)\|_{\g g^*} \le R\},
%$$
%exists for every test form $\phi \in \Omega^{top}_c(\g g)$ and the mapping
%$$
%\int_M^{distrib} \alpha \wedge e^{i\tilde\omega}: \quad
%\Omega^{top}_c(\g g) \to \BB C,
%\quad \phi \mapsto \lim_{R \to +\infty} 
%\int_{M_{\le R}} \Bigl(
%\int_{\g g} \alpha \wedge e^{i\tilde\omega} \wedge \phi \Bigr)
%$$
%is continuous.
Moreover, the distribution $\int_M^{distrib} \alpha \wedge e^{i\tilde\omega}$
is tempered, in particular its Fourier transform makes sense.
%$$
%{\cal F} \biggl( \int_M^{distrib} \alpha \wedge e^{i\tilde\omega} \biggr)
%= \int_M^{distrib} {\cal F}(\alpha \wedge e^{i\tilde\omega})^{distrib}.
%$$
\end{prop}

This proposition follows from Lemma \ref{polynomial1} and a simple
observation that the integral
$$
\int_{\g g} \alpha(X) \wedge e^{i\langle \xi, X \rangle} \wedge \phi(X),
\qquad \xi \in \g g^*,
$$
decays rapidly as $\xi \to \infty$,
for every test form $\phi \in \Omega^{top}_c(\g g)$.

%In his paper \cite{Par2} P.-E.~Paradan studies the distribution
%$\int_M^{distrib} {\cal F}(\alpha \wedge e^{i\tilde\omega})^{distrib}$ on $\g g^*$
%and gives a localization formula for it at (the connected components of) the
%set critical points of $\|\mu\|_{\g g^*}^2$.
%In the subanalytic setting the function $\|\mu\|_{\g g^*}^2$, the set
%$Cr \|\mu\|_{\g g^*}^2$ and its set of critical values are subanalytic.
%By Sard's theorem the set of critical values has measure zero,
%hence by o-minimality it has to be finite.
%Since $\mu$ is proper, the set $Cr \|\mu\|_{\g g^*}^2$ is compact.
%Then $Cr \|\mu\|_{\g g^*}^2$ has finitely many connected components,
%each of which is a subanalytic set.

\separate

\begin{section}
{The fixed point localization formula}  \label{loc1sect}
\end{section}

\subsection
{The first fixed point localization formula}

Let $(M,\omega,\mu)$ be a subanalytic Hamiltonian system in the sense of
Definition \ref{subanalytic-hamiltonian-def}
%From now on we assume that
%``everything is subanalytic,'' i.e. the manifold $M$,
%the symplectic form $\omega$, the moment map $\mu$,
%the action of the group $G$ on $M$ and
and $\alpha \in \Omega_G^*(M)$ a subanalytic equivariant form in the
sense of Definition \ref{subanalytic} and discussion at the end of
Subsection \ref{structure-intro}.
The moment map $\mu$ is always assumed to be proper.

Let $T \subset G$ denote the maximal torus with Lie algebra $\g t$,
and denote by $M^T$ the set of points in $M$ fixed by $T$.
Following \cite{HP} we say that the action of $G$ on $M$ is
{\em $T$-compact} if the set $M^T$ is compact.
Note that since all maximal tori $T \subset G$ are conjugate, if
the set $M^T$ is compact for one particular torus then $M^T$ is compact
for all tori, so the choice of a torus is irrelevant here.
In this subsection we study the localization properties of
$\int_M^{distrib} \alpha \wedge e^{i\tilde\omega}$
at the fixed point set $M^T$.

For each maximal torus $T \subset G$, we denote by ${\mathcal N}_{M^T}$
the normal bundle at $M^T$.
Then we denote by $\chi_T({\mathcal N}_{M^T}) \in \Omega_T^*(M^T)$
the $T$-equivariant Euler form of ${\mathcal N}_{M^T}$, it is a map
$\g t \to \Omega^*(M^T)$.
The $T$-equivariant Euler form is a concrete differential form realization
of the Euler class in $H_T^*(M^T)$ and is determined up to an exact form.
(See, for instance, \cite{BGV} for details.)

We denote by $\g g^{rs}$ the set of {\em regular semisimple} elements
in $\g g$.
These are elements $X \in \g g$ such that the adjoint action of
$ad(X)$ on $\g g$ is diagonalizable (over $\BB C$) and has maximal
possible rank. The set $\g g^{rs}$ is an open and dense subset of $\g g$.

Next we introduce the set of {\em strongly regular elements} $\g g'$.
It consists of regular semisimple elements
$X \in \g g^{rs}$ which satisfy the following additional properties.
If $\g t(X) \subset \g g$ is the unique Cartan
subalgebra in $\g g$ containing $X$, then:
\begin{enumerate}
\item
The set of zeroes $M^X$ is exactly the set of points in $M$ fixed by
the torus $T(X) = \exp (\g t(X)) \subset G$:
$$
M^X = M^{T(X)};
$$
\item
The component of the equivariant Euler form
$$
\chi_{T(X)}({\mathcal N}_{M^X})_{[0]}(X) \ne 0
$$
(i.e. $\chi_{T(X)}({\mathcal N}_{M^X})(X)$ is invertible) at all point $m \in M^X$.
\end{enumerate}
Clearly, $\g g'$ is an open $Ad(G)$-invariant subset of $\g g$.
Since $M^T$ is compact, it has finitely many connected components.
So, for any Cartan subalgebra $\g t \subset \g g$, the intersection
$\g g' \cap \g t$ is just $\g t$ without a finite number of hyperplanes.
Hence, by $Ad(G)$-invariance, the complement of $\g g'$ in $\g g$
has measure zero and $\g g'$ is dense in $\g g$.
When $X \in \g g' \cap \g t$, we have $\g t(X) = \g t$, $T(X)=T$
and $M^X=M^{T(X)}=M^T$.

We need to ensure that, for each $X \in \g g'$, the vector fields $\VF_X$
do not have ``zeroes at infinity.''
It is not enough to say there is some Riemannian metric $(\cdot,\cdot)_M$
on $M$ such that the function $(\VF_X,\VF_X)_M$ is bounded away from zero
on the complement of a compact subset of $M$ containing $M^X$, as 
$(\cdot,\cdot)_M$ can always be scaled.
Thus we need to require $(\VF_X,\VF_X)_M$ to be bounded away from zero
relatively to the coefficients of the metric $(\cdot,\cdot)_M$ itself.
%If $(\cdot,\cdot)_M$ is a metric on $M$, we denote by $g_* (\cdot,\cdot)_M$
%the metric on $M$ obtained by translation by the action of $g \in G$.
%Note that if we start with a subanalytic metric on $M$, its average by the
%$G$-action need not be subanalytic.
Fix a vector space basis $\{E^1,\dots,E^{\dim \g g}\}$ of $\g g$ and an open set
$U \subset M$ containing $G \cdot M^T$ and having compact closure.

\begin{df}  \label{no-zeroes}
We say that the action of $G$ on $M$ {\em has no zeroes at infinity} if it is
$T$-compact and there exists a subanalytic Riemannian metric
$(\cdot,\cdot)_M$ on $M$ with the following property:
For each compact subset $D \subset \g g'$, there is a constant $c_D >0$
such that
\begin{equation}  \label{no-zeroes-eq}
(\VF_X,\VF_X)_M(m) \ge c_D \quad \text{and} \quad
(\VF_X,\VF_X)_M(m) \ge c_D \cdot \bigl| (\VF_{E^a},\VF_{E^b})_M(m) \bigr|,
%\qquad
%\forall X \in D, \quad \forall m \in M \setminus U, \quad \forall g \in G,
\end{equation}
for all $X \in D$, $m \in M \setminus U$ and $1 \le a,b \le \dim \g g$.
\end{df}

%Note that if we start with a subanalytic metric on $M$, its average by the
%$G$-action need not be subanalytic.
The notion of the action of $G$ on $M$ having no zeroes at infinity does not
depend on the choice of vector space basis of $\g g$, nor does it depend
on the choice of open set $U \subset M$ containing $G \cdot M^T$
and having compact closure.

The  metric $(\cdot,\cdot)_M$ in the definition is not required to be
$G$-invariant.
Note that if we start with a subanalytic metric on $M$, its average by the
$G$-action need not be subanalytic.
For a metric $(\cdot,\cdot)_M$ on $M$, we denote by $g_* (\cdot,\cdot)_M$
the metric obtained by translation by the action of $g \in G$.
Without loss of generality we can assume that $U \subset M$ is $G$-invariant.
Then one can choose the constants $c_D$ so that
(\ref{no-zeroes-eq}) will be satisfied with $g_* (\cdot,\cdot)_M$ in place of
$(\cdot,\cdot)_M$, for all $g \in G$.

\begin{thm} \label{loc1}
Let $(M,\omega,\mu)$ be a subanalytic Hamiltonian system,
the manifold $M$ need not be compact.
Suppose that the action of $G$ has no zeroes at infinity, and that the
moment map $\mu: M \to \g g^*$ is proper.
Let $\alpha \in \Omega_G^*(M)$ be a subanalytic
equivariant form which is equivariantly closed.
Then the restriction of the distribution
$\int_M^{distrib} \alpha \wedge e^{i\tilde\omega}$ to $\g g'$
is an $Ad(G)$-invariant function on $\g g'$
\begin{equation}  \label{mainequation1}
F_{\alpha,\tilde\omega}(X) = (-2\pi i)^n
\int_{M^X} \biggl( \frac {(\alpha(X) \wedge e^{i\tilde\omega})|_{M^X}}
{\chi_{T(X)}({\mathcal N}_{M^X})(X)} \biggr)_{[\dim M^X]},
\qquad X \in \g g',
\end{equation}
where $n= \frac 12 \dim M$. That is, if $\phi \in \Omega_c^{top}(\g g')$
is a smooth compactly supported differential form on $\g g'$ of top degree,
$$
\int_M^{distrib} \alpha \wedge e^{i\tilde\omega}: \quad
\phi \mapsto \int_{\g g} F_{\alpha,\tilde\omega} \phi.
$$
%We extend the function $F_{\alpha,\tilde\omega}$ by zero to a measurable
%function on $\g g$. If $F_{\alpha,\tilde\omega}$ happens to be locally
%integrable with respect to the Lebesgue measure on
%$U \subset \g g \simeq \BB R^{\dim \g g}$,
%then the equation (\ref{mainequation}) holds for smooth differential
%forms $\phi$ of top degree which are compactly supported on $U$
%(and not necessarily on $U \cap \g g'$).
\end{thm}

Distribution $\int_M^{distrib} \alpha \wedge e^{i\tilde\omega}$ were
considered by the author in \cite{L1,L2}.
However, the papers \cite{L1, L2} deal with actions of non-compact groups
preserving a homology cycle in $T^*M$, so the group and the manifold are much
more general, but the Hamiltonian structure is the one of the cotangent bundle.
In the special case when $M$ is a coadjoint orbit, Theorem \ref{loc1}
was proved by M.~Duflo, G.~Heckman and M.~Vergne \cite{DHV}, \cite{DV}
and later by P.-E.~Paradan \cite{Par1}.

\begin{rem}
Returning to Example \ref{ex1} where the equivariant form $\alpha=1$,
this result formally coincides with the Duistermaat-Heckman formula \cite{DH}.
If, in addition, the manifold $M$ is compact, this result is exactly
the Duistermaat-Heckman formula.
\end{rem}

This localization formula (\ref{mainequation1}) is also closely related
to the extension of the Duistermaat-Heckman formula to non-compact manifolds
due to E.~Prato and S.~Wu \cite{PrWu}.
They consider the integral $\int_M^{distrib} e^{i\tilde \omega}$
(i.e. $\alpha =1$), but they {\em do not} work in the subanalytic setting.
To work around the problem of convergence of this integral they assume
that there exists an $X_0 \in \g g'$ such that the component
of the moment map
$$
\mu_{X_0}: M \to \BB R, \qquad \mu_{X_0} =_{def} \langle \mu, X_0 \rangle
$$
is proper and not surjective. This implies that $\mu_{X_0}$ is
{\em polarized}, i.e. bounded either from below or from above.
This assumption is similar to our significantly weaker assumption that the
moment map $\mu: M \to \g g^*$ is proper.
Like us, they require the group action to be $T$-compact
(although they do not use this term).
Finally, they complexify the Lie algebra $\g t(X_0)$ and prove that, for
$Z \in \g t(X_0) \otimes_{\BB R} \BB C$ with $\re (Z) \in \g t(X_0) \cap \g g'$
and $\im (Z)$ lying in the interior of a certain cone
${\cal C} \subset \g t(X_0)$, the improper integral
$\int_M e^{i(\mu(Z)+\omega)}$ converges to
$$
F_{1,\tilde\omega}(Z) =
\int_{M^{T(X_0)}} \biggl( \frac {e^{i(\mu(Z)+\omega)}|_{M^{T(X_0)}}}
{\chi_{T(X_0)}({\mathcal N}_{M^{T(X_0)}})(Z)} \biggr)_{[\dim M^{T(X_0)}]}
$$
in the most common sense of convergence (and in particular in the sense
of distributions).
Note that $0 \in \partial {\cal C}$ and the interior of ${\cal C}$ being
non-empty is essentially equivalent to $\mu_{X_0}$ being polarized.

\subsection
{Proof of Theorem \ref{loc1}}

If we knew in addition that the moment map $\mu : M \to \g g^*$ composed with
the projection $\g g^* \twoheadrightarrow \g t^*$ was proper and that all
integrals converged, then the classical argument of
N.~Berline and M.~Vergne \cite{BGV} would apply verbatim.
However, we only assume that the moment map $\mu$ itself is proper,
so we cannot deal with integral (\ref{mainequation1})
``one Cartan algebra at a time'' and we proceed as in \cite{BV2}.

Fix a subanalytic metric $(\cdot,\cdot)_M$ on $M$ satisfying conditions of
Definition \ref{no-zeroes}, and let $(\cdot,\cdot)_M^G$ be its average by
$G$-action.
Then we define a 1-form $\theta_X$ depending on $X \in \g g$ by setting
$$
\theta_X = \frac {( \VF_X,\,\cdot\,)_M^G}{(\VF_X,\VF_X)_M^G}.
$$
%$$
%\theta_X = \frac {( \VF_X,\,\cdot\,)_M}{(\VF_X,\VF_X)_M}.
%$$
For a fixed $X \in \g g$, this form is defined on $M \setminus M^X$.
We regard $\theta_X$ as a map
$\{(X,m) \in \g g \times M ;\: m \notin M^X\} \to \Lambda^*(TM)$.
Note that
$$
\{(X,m) \in \g g \times M ;\: m \notin M^X\} \quad \supset \quad
\g g' \times (M \setminus G \cdot M^T),
$$ 
where $T \subset G$ is a maximal torus, and the set $G \cdot M^T$ is compact.
The form $\theta_X$ is $G$-equivariant and has the following property
$$
\iota_X \theta_X =1, \qquad \forall X \in \g g.
$$
Hence, on $\{(X,m) \in \g g \times M ;\: m \notin M^X\}$, we have
$(d_{\g g} \theta_X)_{[0]}=1$, the form $d_{\g g} \theta_X$ is invertible
and the quotient $\frac{\theta_X}{d_{\g g}\theta_X}$ makes sense.

Recall that the test form $\phi$ is compactly supported in $\g g'$.
For each $R > 0$ large enough so that $G \cdot M^T$ is contained in the
interior of $M_{\le R}$, by the classical localization argument we have:
\begin{multline*}
\int_{M_{\le R}} \Biggl(
\int_{\g g} \alpha \wedge e^{i\tilde\omega} \wedge \phi \Biggr)
= \int_{\g g} \Biggl(
\int_{M_{\le R}} \alpha \wedge e^{i\tilde\omega} \Biggr) \wedge \phi  \\
=\int_{\g g} \int_{M^X} \Biggl( \frac{(\alpha \wedge e^{i\tilde\omega})|_{M^X}}
{\chi_{T(X)}({\mathcal N}_{M^X})(X)} \biggr)\wedge \phi+
\int_{\g g} \biggl( \int_{\partial M_{\le R}}
\frac{\theta_X}{d_{\g g}\theta_X} \wedge \alpha \wedge e^{i\tilde\omega}
\biggr) \wedge \phi  \\
= \int_{\g g} F_{\alpha,\tilde\omega} \phi +
\int_{\partial M_{\le R}} \biggl( \int_{\g g}
\frac{\theta_X}{d_{\g g}\theta_X} \wedge \alpha \wedge e^{i\tilde\omega}
\wedge \phi \biggr)_{[2n-1]}.
\end{multline*}
As $R \to +\infty$, the left hand side tends to
$\bigl\langle \int_M^{distrib} \alpha \wedge e^{i\tilde\omega}, \phi \bigr\rangle$,
so it remains to show
\begin{equation}  \label{lim=0}
\int_{\partial M_{\le R}} \biggl( \int_{\g g}
\frac{\theta_X}{d_{\g g}\theta_X} \wedge \alpha \wedge e^{i\tilde\omega}
\wedge \phi \biggr)_{[2n-1]}
\quad \to 0 \qquad \text{as $R \to +\infty$}.
\end{equation}

Let $h: \g g \times M \to \BB R$ be the function $h(X,m)= (\VF_X,\VF_X)_M^G(m)$.
The form $\alpha(X) \in \Omega_G^*(M)$ depends on $X \in \g g$ polynomially.
Hence the form
$\frac{\theta_X}{d_{\g g}\theta_X} \wedge \alpha(X) \wedge e^{i\omega}$
depends on the parameter $X \in \g g$ as a polynomial in $X$ and $1/h(X,m)$,
thus it can be expressed as
$$
\frac{\theta_X}{d_{\g g}\theta_X} \wedge \alpha(X) \wedge e^{i\omega} =
\sum_{j=1}^k \beta_j \cdot \frac{P_j(X)}{h(X,m)^{n_j}}
$$
for some polynomials $P_j(X)$ on $\g g$, $n_j \ge 0$ and
$\beta_j \in \Omega^*(M)$ obtained by averaging subanalytic forms by
$K$-actions.
By Lemma \ref{polynomial2}, there is an $N \in \BB N$ such that each
$$
\int_{\partial M_{\le R}} |\beta_j| \quad
\text{is $O(R^N)$ as $R \to +\infty$}, \qquad 1 \le j \le k.
$$
With $D=\supp \phi$ in Definition \ref{no-zeroes}, condition
(\ref{no-zeroes-eq}) implies that the partial derivatives of all orders
of $1/h(X,m)^{n_j}$ with respect to $X \in \g g$ are uniformly bounded.
More precisely, for each multiindex $L=(l_1, \dots,l_{\dim \g g})$,
there exists a constant $C_L>0$ such that
$$
%\frac{\partial^{|L|}}{(\partial E^1)^{l_1} \dots (\partial E^{\dim \g g})^{l_{\dim \g g}}}
\frac{\partial^{|L|}}{\partial E^L} \biggl( \frac 1{h(X,m)^{n_j}} \biggr) \le C_L,
\qquad \forall X \in \supp \phi, \quad \forall m \in M \setminus U.
$$
%for all $X \in \supp \phi$ and $m \in M \setminus U$.
Hence the Fourier transforms
$$
\int_{\g g} \frac{P_j(X)}{h(X,m)^{n_j}} \cdot \phi(X) \cdot e^{i\xi(X)},
\qquad \xi \in \g g^*, \qquad 1 \le j \le k,
$$
are $o(R^{-N})$ as $\|\xi\|_{\g g^*}=R \to +\infty$.
This proves (\ref{lim=0}).

It is clear that the function $F_{\alpha,\tilde\omega}$ is $Ad(G)$-invariant.
\qed

\subsection
{Definition of $\int_M^{distrib} \alpha$ and the main localization theorem}
\label{regularization}

Note that, for each $s \in \BB R$, $s \ne 0$, the pair $(s\omega,s\mu)$
gives another symplectic structure on $M$ such that the action of $G$
remains Hamiltonian.
As before, consider a subanalytic equivariant form
$\alpha \in \Omega_G^*(M)$ which is equivariantly closed.
Then Theorem \ref{loc1} applied to $(M,s\omega,s\mu)$ implies that,
for each $\phi \in \Omega^{top}_c(\g g')$, the limit
$$
\lim_{s \to 0^+}
\biggl\langle \int_M^{distrib} \alpha \wedge e^{is\tilde\omega},
\phi \biggr\rangle
$$
exists and the assignment
$$
\phi \mapsto \lim_{s \to 0^+}
\biggl\langle \int_M^{distrib} \alpha \wedge e^{is\tilde\omega},
\phi \biggr\rangle
$$
is a distribution on $\g g'$.
Recall that the space of distributions is equipped with the weak*-topology,
and a sequence of distributions $\{\Lambda_j\}$ converges to a
distribution $\Lambda$ in this topology if and only if
$\lim_{j \to \infty} \Lambda_j (\phi) = \Lambda (\phi)$ for all test
functions $\phi$.
Therefore, we can define $\int_M^{distrib} \alpha$
as a limit of distributions on $\g g'$ (in the weak*-topology):
$$
\int_M^{distrib} \alpha = \lim_{s \to 0^+}
\int_M^{distrib} \alpha \wedge e^{is\tilde\omega}.
$$

\begin{rem}
We do not allow $s<0$ in the limit because the orientation of $M$
is determined by its symplectic structure and replacing $\omega$ with
$s\omega$, $s<0$, will change the orientation whenever
$n = \frac 12 \dim M$ is odd.
\end{rem}

The following localization formula follows immediately from Theorem \ref{loc1}:

\begin{thm} \label{main}
Let $(M,\omega,\mu)$ be a subanalytic Hamiltonian system, the manifold $M$
need not be compact.
Suppose that the action of $G$ has no zeroes at infinity and the
moment map $\mu: M \to \g g^*$ is proper.
Let $\alpha \in \Omega_G^*(M)$ be a subanalytic
equivariant form which is equivariantly closed.
Then the distribution $\int_M^{distrib} \alpha$ on $\g g'$
is given by integrating against a function $F_{\alpha}$:
$$
\int_M^{distrib} \alpha : \quad \phi \mapsto \int_{\g g} F_{\alpha} \phi,
\qquad \qquad \phi \in \Omega^{top}_c(\g g')
$$
where $F_{\alpha}$ is an $Ad(G)$-invariant function on $\g g'$
given by the formula
\begin{equation*}
F_{\alpha}(X) = (-2\pi i)^n
\int_{M^X} \biggl( \frac {\alpha(X)|_{M^X}}
{\chi_{T(X)}({\mathcal N}_{M^X})(X)} \biggr)_{[\dim M^X]},
\qquad X \in \g g',
\end{equation*}
where $n= \frac 12 \dim M$.
\end{thm}

Note that this result formally coincides with the classical
Berline-Vergne localization formula \cite{BV1}, \cite{BGV}.

\begin{cor}  \label{cor1}
If $\alpha \in \Omega_G^*(M)$ is equivariantly exact, i.e.
$\alpha = d_{\g g} \beta$ for some $\beta \in \Omega_G^*(M)$, then
$$
\int_M^{distrib} \alpha = 0
$$
as a distribution on $\g g'$.
Hence the map $\alpha \mapsto \int_M^{distrib} \alpha \bigr|_{\g g'}$
descends to cohomology $H_G^*(M)$.
\end{cor}

\pf
By the localization formula (Theorem \ref{main})
it is sufficient to prove that the function $F_{\alpha}=0$ on $\g g'$.
Note that for each $X \in \g g'$, the vector field $\VF_X$ is zero on $M^X$, so
\begin{multline*}
(-2\pi i)^{-n} \cdot F_{\alpha}(X) =
\int_{M^X} \biggl( \frac {\alpha(X)|_{M^X}}
{\chi_{T(X)}({\mathcal N}_{M^X})(X)} \biggr)_{[\dim M^X]}
= \int_{M^X} d_{\g g} \biggl( \frac {\beta(X)|_{M^X}}
{\chi_{T(X)}({\mathcal N}_{M^X})(X)} \biggr)  \\
= \int_{M^X} d \biggl( \frac {\beta(X)|_{M^X}}
{\chi_{T(X)}({\mathcal N}_{M^X})(X)} \biggr)_{[\dim M^X -1]}
= \int_{\partial M^X} \biggl( \frac {\beta(X)|_{M^X}}
{\chi_{T(X)}({\mathcal N}_{M^X})(X)} \biggr)_{[\dim M^X -1]} =0.
\end{multline*}
\qed

\begin{rem}
Theorem \ref{main} together with Corollary \ref{cor1} provide an alternative
approach to the definition of an integral given by
T.~Hausel and N.~Proudfoot in \cite{HP}.
They consider the equivariant cohomology $H^*_G(M)$,
tensor it with the field of $Ad(G)$-invariant rational functions on $\g g$:
$$
\widehat{H}^*_G(M) = H^*_G(M) \otimes_{\BB R [\g g]^G} \BB R (\g g)^G,
$$
and they want to make sense out of the integral $\int_M \alpha$,
where $\alpha \in \widehat{H}^*_G(M)$ and the manifold $M$ is not compact.
So under an additional assumption that the group action is $T$-compact
they {\em define}
$$
\int_M \alpha(X) = (-2\pi i)^n
\int_{M^X} \biggl( \frac {\alpha(X)|_{M^X}}
{\chi_{T(X)}({\mathcal N}_{M^X})(X)} \biggr)_{[\dim M^X]},
\qquad X \in \g g.
$$
%This is a valid definition and it does not need any further justification.
From our point of view, both sides exist and equal as distributions on a
dense open $Ad(G)$-invariant subset $\g g'$ of $\g g$.
%(Note that $\Omega(M,G)$ is a subcomplex of the standard twisted deRham
%complex; $\Omega(M,G)$ consists of those differential forms
%$\alpha(X)$, $X \in \g g$, which are real algebraic on $M$.
%This subcomplex also computes the equivariant cohomology $H^*_G(M)$.)
\end{rem}

%We conclude this subsection with a discussion related to P.-E.~Paradan's
%paper \cite{Par2}. We continue to assume that $(M,\omega,\mu)$ is a
%subanalytic Hamiltonian system. Then the function $\|\mu\|_{\g g^*}^2$,
%its set of critical points $Cr \|\mu\|_{\g g^*}^2$ and its set of critical
%values are subanalytic.
%By Sard's theorem the set of critical values has measure zero,
%hence by o-minimality it has to be finite.
%Since $\mu$ is proper, the set $Cr \|\mu\|_{\g g^*}^2$ is compact.
%Then $Cr \|\mu\|_{\g g^*}^2$ has finitely many connected components,
%each of which is a subanalytic set.

%Let ${\cal H}$ be the Hamiltonian vector field on $M$ associated to
%$\frac 12 \|\mu\|_{\g g^*}^2$. (This means that ${\cal H}$ is the unique vector
%field on $M$ satisfying $\frac 12 d(\|\mu\|^2) = \iota({\cal H}) \omega$.)
%Fix a $G$-invariant Riemannian metric $(\cdot,\cdot)_M$ on $M$.
%The Witten 1-form is a $G$-invariant form on $M$ defined by
%$$
%\lambda =_{def} ({\cal H}, \cdot )_M.
%$$
%Denote by $\Phi_{\lambda}: M \to \mathfrak{g}^*$ the $G$-invariant
%map defined by
%$$
%\langle \Phi_{\lambda}(m), X \rangle =_{def} \lambda_m (\VF_X \bigr|_m),
%\qquad m \in M, \quad X \in \g g.
%$$
%Then $\{ \Phi_{\lambda} = 0 \}$ is precisely the set of critical points
%of $\|\mu\|^2$:
%$$
%\{ m \in M;\: \Phi_{\lambda}(m) = 0 \} = \operatorname{Cr} (\|\mu\|^2).
%$$

\subsection
{Coadjoint orbits}

In this subsection we illustrate how our integration theory for non-compact
subanalytic manifolds applies to coadjoint orbits of real semisimple
Lie groups. This special case is very important and it has been thoroughly
studied by many mathematicians.
Just to mention a few works on this subject, see \cite{BV2}, \cite{DHV},
\cite{DV}, \cite{Par1}, \cite{Ro}, \cite{S} and references therein.

Let $G^{ss}$ be a real semisimple Lie group and $G \subset G^{ss}$
its maximal compact subgroup. We denote by $\g g^{ss}$ and $\g g$
the Lie algebras of $G^{ss}$ and $G$ respectively.
Let $\lambda \in (\g g^{ss})^*$ be a semisimple element,
and consider its coadjoint orbit
$$
{\cal O}_{\lambda} = G^{ss} \cdot \lambda \subset (\g g^{ss})^*.
$$
When $\lambda \in (\g g^{ss})^*$ is semisimple, the orbit ${\cal O}_{\lambda}$
is a closed submanifold of $(\g g^{ss})^*$.
%Recall a well-known result from symplectic geometry.
%\begin{prop}
%The coadjoint orbit ${\cal O}_{\lambda}$ has a canonical $G^{ss}$-invariant
%symplectic form given by
%$$
%\omega (\VF_X,\VF_Y)_f = -f([X,Y]), \qquad f \in {\cal O}_{\lambda}, \quad
%X, Y \in \g g^{ss}.
%$$
%Furthermore, the action of $G^{ss}$ on ${\cal O}_{\lambda}$ is Hamiltonian,
%with symplectic moment map given by the inclusion map
%${\cal O}_{\lambda} \hookrightarrow (\g g^{ss})^*$.
%\end{prop}
Recall that the coadjoint orbit ${\cal O}_{\lambda}$ has the Konstant-Kirillov
symplectic form which is $G^{ss}$-invariant, and the action of $G^{ss}$ on
${\cal O}_{\lambda}$ is Hamiltonian, with symplectic moment map given by the
inclusion map ${\cal O}_{\lambda} \hookrightarrow (\g g^{ss})^*$.
The $G^{ss}$ action on ${\cal O}_{\lambda}$ restricts to a Hamiltonian
action of $G$, with symplectic moment map $\mu: {\cal O}_{\lambda} \to \g g^*$
given by the inclusion ${\cal O}_{\lambda} \hookrightarrow (\g g^{ss})^*$
composed with the natural projection $(\g g^{ss})^* \twoheadrightarrow \g g^*$.
The coadjoint orbit ${\cal O}_{\lambda}$ is a smooth real affine variety.
In order to apply our integration results (Theorems \ref{loc1} and
\ref{main}) we need to know if the moment map $\mu$ is proper and
the $G$-action is $T$-compact.

\begin{prop}
Let $G^{ss}$ be a real semisimple Lie group and $G \subset G^{ss}$
a maximal compact subgroup.
Denote by $\g g^{ss}$ and $\g g$ their respective Lie algebras.
Let ${\cal O}_\lambda \subset (\g g^{ss})^*$ be a semisimple
coadjoint orbit of $G^{ss}$.
Then the restriction of the natural projection map
$(\g g^{ss})^* \twoheadrightarrow \g g^*$ to ${\cal O}_{\lambda}$ is proper.

Let $T \subset G$ be a maximal torus with Lie algebra $\g t$.
If $\g t$ is also a Cartan algebra in $\g g^{ss}$, then $({\cal O}_{\lambda})^T$
-- the set of points in ${\cal O}_{\lambda}$ fixed by $T$ -- is compact.
In fact, either ${\cal O}_{\lambda}$ is an elliptic orbit
(i.e. ${\cal O}_{\lambda} \cap \g t^* \ne \varnothing$)
and $({\cal O}_{\lambda})^T$ is finite or ${\cal O}_{\lambda}$ is not elliptic
and $({\cal O}_{\lambda})^T$ is empty.
\end{prop}

This is a well-known result and its proofs can be found, for instance,
in \cite{DHV}, \cite{L3}.
In the setting of closed coadjoint orbits and $\alpha =1$,
Theorem \ref{loc1} essentially reduces to Kirillov's character formula
due to W.~Rossmann \cite{Ro},
and a proof similar to ours has originally appeared in \cite{BV2}.
Note that in this case Theorem \ref{loc1} describes the restriction of the
distribution $\int_M^{distrib} \alpha \wedge e^{i\tilde\omega}$ to $\g g^{rs}$.
In the case of coadjoint orbits one can compute the entire (unrestricted)
distribution $\int_M^{distrib} \alpha \wedge e^{i\tilde\omega}$ on $\g g$,
as was done by P.-E.~Paradan in \cite{Par1}.
(Earlier work in this direction was done by J.~Sengupta \cite{S},
and by M.~Duflo, G.~Heckman and M.~Vergne \cite{DHV} and \cite{DV}.)

One might ask if the natural projection $\g g^* \twoheadrightarrow \g t^*$
is proper on ${\cal O}_{\lambda}$.
E.~Prato (Propositions 2.2 and 2.3 in \cite{Pr}) shows that when
$(G^{ss}, G)$ is an irreducible Hermitian symmetric pair, the coadjoint orbit
${\cal O}_{\lambda}$ is elliptic and $\lambda$ lies in a certain cone of
$\g t^* \subset (\g g^{ss})^*$ the answer to the question is affirmative.
More precisely, it is known that when $(G^{ss}, G)$ is an irreducible Hermitian
symmetric pair, the center $Z \subset G$ is the circle group \cite{He},
and E.~Prato shows that there exists an element $X_0 \in Lie(Z)$
such that the component of the moment map
$$
\mu_{X_0}: M \to \BB R, \qquad \mu_{X_0} = \langle \mu, X_0 \rangle
$$
is proper and bounded from below, hence polarized.
But in general the projection $(\g g^{ss})^* \twoheadrightarrow \g t^*$
need not be proper on ${\cal O}_{\lambda}$.
See \cite{L3} for a concrete counterexample.

\separate

\begin{section}
{Appendix: O-minimal structures and subanalytic sets}
\end{section}

Examples of subanalytic sets include complex affine varieties and sets in
$\BB R^n$ defined by finitely many real polynomial equations and inequalities.
The collection of all subanalytic sets is a particular example of a structure.
In this section we give a very brief introduction to structures -- collections
of sets in $\BB R^n$ with nice geometric properties.
Structures are studied in model theory, a very exciting part of logic which
produces extremely useful and highly non-trivial geometric results.
For more details, proofs and further references the reader is referred to
the works by L. van den Dries and C. Miller \cite{DM} and \cite{D}.
In this section we summarize \cite{DM} and list the key properties of
o-minimal structures which imply Lemmas \ref{polynomial1}, \ref{polynomial2}
which in turn are used to prove Proposition \ref{subanalytic-hamiltonian-prop}
and Theorem \ref{loc1}.

\subsection
{Definitions and basic properties}  \label{structure-intro}

\begin{df}  \label{str}
A structure (on the real field $(\BB R,+,\cdot)$) is a sequence
$\mathfrak S = \{ \mathfrak S_n \}_{n \in \BB N}$
such that for each $n \in \BB N$:
\begin{enumerate}
\item
$\mathfrak S_n$ is a boolean algebra of subsets of
$\BB R^n$, with $\BB R^n \in \mathfrak S_n$;
\item 
$\mathfrak S_n$ contains the diagonal
$\{ (x_1,\dots,x_n) \in \BB R^n ;\: x_i = x_j\}$ for $1 \le i < j \le n$;
\item
If $A \in \mathfrak S_n$, then $A \times \BB R$ and $\BB R \times A$ belong
to $\mathfrak S_{n+1}$;
\item {\label{4}}
If $A \in \mathfrak S_{n+1}$, then $\pi(A) \in \mathfrak S_n$, where
$\pi: \BB R^{n+1} \to \BB R^n$ is the projection on the first $n$ coordinates;
\item
$\mathfrak S_3$ contains the graphs of addition and multiplication.
\end{enumerate}
\end{df}

We say that a set $A \subset \BB R^m$ belongs to $\mathfrak S$
(or $\mathfrak S$ contains $A$) if $A \in \mathfrak S_m$; and that
a (not necessarily continuous) map $f: A \to \BB R^n$ belongs to
$\mathfrak S$ (or $\mathfrak S$ contains $f$) if the graph of $f$ lies in
$\mathfrak S_{m+n}$.

Although part \ref{4} of the definition is not symmetric with respect to the
coordinates $x_1,\dots,x_n$, it follows that $\mathfrak S$ is invariant
under ``permutations, repetitions and omission of the coordinates''.
That is, if $B \in \mathfrak S_n$ and $i_1,\dots,i_n \in \{1,\dots,m\}$
(repetitions allowed), then the set
$$
\{ (x_1,\dots,x_m) \in \BB R^m ;\: (x_{i_1},\dots,x_{i_n}) \in B \}
$$
belongs to $\mathfrak S_m$.

For each structure $\mathfrak S$, $\mathfrak S_1$ automatically contains all
singleton sets $\{a\}$, $a \in \BB Q$, and open and closed intervals
with rational endpoints, $\mathfrak S_2$ contains the set
$\{(x,y)\in \BB R^2 ;\: x<y \}$ which may be interpreted as
``the order relation $<$ belongs to $\mathfrak S$''.
For each polynomial with rational coefficients
$f(X_1,\dots,X_n) \in \BB Q[X_1,\dots,X_n]$,
the corresponding function $f: \BB R^n \to \BB R$,
$x \mapsto f(x)$, belongs to $\mathfrak S$.
If $A \in \mathfrak S_m$, then the interior, closure and boundary of $A$
also belong to $\mathfrak S_m$.
If $f: A \to \BB R$ belongs to $\mathfrak S$, then the sets
$\{x \in A ;\: f(x)=0 \}$ and $\{x \in A ;\: f(x) > 0 \}$
also belong to $\mathfrak S$.

Let $A \subset \BB R^m$.
For each $k,p \in \BB N$, let $\operatorname{Reg}_k^p(A)$ denote the set
of all $x \in A$ having an open neighborhood $U$ of $x$ such that
$A \cap U$  is an (embedded) ${\cal C}^p$ submanifold of $\BB R^m$ of
dimension $k$. If $A$ belongs to $\mathfrak S$, then so does each
$\operatorname{Reg}_k^p(A)$.

If $A \subset \BB R^m$, $B \subset \BB R^n$ and the functions
$f: A \to \BB R^n$, $g: B \to \BB R^q$ belong to $\mathfrak S$,
then the composition $g \circ f : f^{-1}(B) \to \BB R^q$ belongs to
$\mathfrak S$.
If $f: A \to \BB R^n$ belongs to $\mathfrak S$ and is injective, then its
compositional inverse $f^{-1}: f(A) \to \BB R^m$ belongs to $\mathfrak S$.
If $f: A \to \BB R^n$ belongs to $\mathfrak S$ and $A$ is open,
then the set of points in $A$ where $f$ is continuous and
the set of points in $A$ where $f$ is differentiable belong to
$\mathfrak S$. If $f$ is differentiable on $A$, then each partial derivative
of $f$ also belongs to $\mathfrak S$.

Suppose that $A$ is a ${\cal C}^1$ submanifold of $\BB R^n$. We identify
both the tangent bundle $T \BB R^n$ and the cotangent bundle $T^* \BB R^n$
with $\BB R^{2n}$ in the obvious way.
In particular, the point $(a,b) \in \BB R^{2n} = T^* \BB R^n$ corresponds
to the linear form $x \mapsto b \cdot x$ on $\BB R^n = T_a \BB R^n$.
These identifications make the tangent bundle $TA$ and the conormal
bundle $T^*_A \BB R^n$ subsets of $\BB R^{2n}$.
Similarly, the exterior bundle $\Lambda^* (T \BB R^n)$ on $T\BB R^n$ can be
identified with $\BB R^{n+2^n}$. This identification makes the exterior
bundle $\Lambda^*(TA)$ a subset of $\BB R^{n+2^n}$.
If $A$ belongs to $\mathfrak S$, then its tangent bundle $TA$, conormal bundle
$T^*_A \BB R^n$ and exterior bundle $\Lambda^*(TA)$ belong to $\mathfrak S$.
In particular, it makes sense to talk about differential forms on $A$
which belong to $\mathfrak S$  -- these are sections
$\alpha: A \to \Lambda^*(TA) \subset \BB R^{n+2^n}$ whose graphs lie in
$\mathfrak S_{2n+2^n}$.

\subsection
{Examples of structures and subanalytic sets}

Given two structures $\mathfrak S$ and $\mathfrak S'$ we put
$\mathfrak S \preccurlyeq \mathfrak S'$ if
$\mathfrak S_n \subset \mathfrak S'_n$ for all $n \in \BB N$; this defines a
partial order on the set of all structures on $(\BB R, +, \cdot)$.

The most trivial (and the least interesting) example of a structure
$\mathfrak S^{max}$ is obtained by letting $\mathfrak S^{max}_n$
be the collection of all subsets of $\BB R^n$.
This is the largest structure on $(\BB R, +, \cdot)$.
We will not consider this structure in this paper.

We denote the smallest structure on $(\BB R,+,\cdot)$ by
$\mathfrak S^{min}$.
Because of the basic properties of structures stated above,
$\mathfrak S^{min}_n$ must contain all finite unions of sets of the form
$$
\{ x \in \BB R^n ;\: f(x)=0, g_1(x)>0, \dots, g_k(x)>0 \}
$$
with $f, g_1,\dots, g_k \in \BB Q[X_1,\dots,X_n]$.
The collection of these finite unions (for $n \in \BB N$) clearly
satisfies conditions 1, 2, 3 and 5 of Definition \ref{str}, and by Tarski's
theorem, also condition 4. Hence $\mathfrak S^{min}$ consists exactly of
these finite unions. A singleton set $\{r\}$ with $r \in \BB R$ belongs to
$\mathfrak S^{min}$ if and only if $r$ is algebraic.

Another example of a structure on $(\BB R,+,\cdot)$ is the collection of
semialgebraic sets denoted by $\mathfrak S^{alg}$. By definition
$\mathfrak S^{alg}_n$ consists of all finite unions of sets of the form
$$
\{ x \in \BB R^n ;\: f(x)=0, g_1(x)>0, \dots, g_k(x)>0 \}
$$
with $f, g_1,\dots, g_k \in \BB R[X_1,\dots,X_n]$.
Like in the previous example, it is clear that $\mathfrak S^{alg}$
satisfies conditions 1, 2, 3 and 5 of Definition \ref{str},
and by Tarski-Seidenberg theorem, also condition 4.

One way to form new structures on $(\BB R,+,\cdot)$ is to pick a
collection of functions $f_j: \BB R^{n_j} \to \BB R$ for $j$ ranging
over some index set $J$ and to consider the smallest structure containing
the graphs of all functions $f_j$. Such structure is called
{\em the structure on $(\BB R, +, \cdot)$ generated by the $f_j$'s}.

\begin{df}
We call a function $f: \BB R^n \to \BB R$ a
{\em restricted analytic function}, if it vanishes identically away from
$[-1,1]^n$ and the restriction of $f$ to $[-1,1]^n$ is analytic.
\end{df}

\begin{df}  \label{subanalytic}
We denote by $\mathfrak S^{an}$ the structure on $(\BB R,+,\cdot)$ generated
by all restricted analytic functions and the functions $x^r : \BB R \to \BB R$
given by
\begin{equation}  \label{x^r}
a \mapsto
\begin{cases}
a^r  & \text{if $a>0$;}  \\
0 & \text{if $a \le 0$,}
\end{cases}
\qquad r \in \BB R.
\end{equation}

We call a set $A \subset \BB R^n$ (respectively a function $f: A \to \BB R^n$)
{\em subanalytic} if the set $A$ (respectively the graph of $f$) belongs to
$\mathfrak S^{an}$.
\end{df}

\begin{rem}
The term ``subanalytic'' is more commonly used to denote sets which belong
to the smaller structure generated by the restricted analytic functions only.
However, for the purposes of this article $\mathfrak S^{an}$ works
just as well and is slightly more general.
\end{rem}

Let $\mathfrak S^{an,\, exp}$ be the structure on $(\BB R,+,\cdot)$ generated
by all restricted analytic functions and $\exp : \BB R \to \BB R$ given by
$\exp(x)=e^x$. Then $\mathfrak S^{an,\, exp}$ contains the logarithm function
$\log : (0,\infty) \to \BB R$, as well as each function $x^r$ defined by
(\ref{x^r}), since $a^r = \exp(r \log a)$ for $a>0$.

We have the following inclusions of structures:
$$
\mathfrak S^{min} \preccurlyeq \mathfrak S^{alg}
\preccurlyeq \mathfrak S^{an} \preccurlyeq \mathfrak S^{an,\, exp}
\preccurlyeq \mathfrak S^{max},
$$
each of these inclusions is strict.
L. van den Dries and C. Miller conjecture \cite{DM} that there are no
structures on $(\BB R,+,\cdot)$ lying strictly between $\mathfrak S^{an}$ and
$\mathfrak S^{an,\, exp}$.

\subsection
{O-minimal structures and their properties}

\begin{df}
A structure $\mathfrak S$ on $(\BB R,+,\cdot)$ is called o-minimal if
$\mathfrak S_1$ consists exactly of the finite unions of intervals of
all kinds (including infinite intervals and singletons).
\end{df}

\begin{ex}
Structures $\mathfrak S^{alg}$, $\mathfrak S^{an}$ and
$\mathfrak S^{an,\, exp}$ are o-minimal, while $\mathfrak S^{min}$
and $\mathfrak S^{max}$ are not.
\end{ex}

O-minimal structures possess particularly nice properties.
We list some of them below in order to demonstrate why it is always preferable
to deal with sets and functions belonging to some o-minimal structure.
From now on we assume that the structure $\mathfrak S$ is o-minimal.

\noindent
{\bf Component theorem.}
Every $A$ belonging to $\mathfrak S$ has finitely many
connected components, each belonging to $\mathfrak S$.
Every connected component of $A$ is also path connected.

\noindent
{\bf Dimension is well-behaved.}
Let $A \in \mathfrak S_n$ be non-empty.
We denote by $\dim A$ the maximum integer $d$ such that $A$ contains a
$d$-dimensional ${\cal C}^1$ submanifold of $\BB R^n$
(so $0 \le \dim A \le n$). We also put $\dim \varnothing = -\infty$.
Then:
\begin{itemize}
\item
$\dim (\overline{A} \setminus A) < \dim A$, where $\overline{A}$
denotes the closure of $A$;
\item
If $f: A \to \BB R^m$ belongs to $\mathfrak S$, then $\dim f(A) \le \dim A$.
\end{itemize}

\noindent
{\bf Monotonicity theorem.}
Let $f : (a,b) \to \BB R$ belong to $\mathfrak S$,
$-\infty \le a < b \le \infty$, and $p \in \BB N$.
Then there are $a_0, a_1,\dots,a_{k+1}$ with
$a=a_0 < a_1 < \dots < a_k < a_{k+1} = b$ such that the restriction of $f$
to each interval $(a_i, a_{i+1})$ is ${\cal C}^p$ and either constant or
strictly monotone, for $i =0,\dots,k$.

\noindent
{\bf Differentiability.}
Let $f : (a,b) \to \BB R$ belong to $\mathfrak S$,
$-\infty \le a < b \le \infty$,
then $f$ is differentiable at all but finitely many points of $(a,b)$.

\noindent
{\bf Triangulation.}
Let $A,A_1,\dots,A_l \in \mathfrak S_n$ with $A_1,\dots,A_l \subset A$.
Then there exist a finite simplicial complex $K$ in $\BB R^n$ and a map
$\phi: A \to \BB R^n$ belonging to $\mathfrak S$ such that $\phi$ maps
$A$  and each $A_i$ homeomorphically onto a union of open simplices of $K$.

\noindent
{\bf Uniform bounds on growth.}
Let $ A \subset \BB R^n$ and $g: A \times \BB R \to \BB R$ belong to
$\mathfrak S$. Then there exist functions $\psi: \BB R \to \BB R$ and
$\rho: A \to \BB R$ belonging to $\mathfrak S$ such that
$|g(x,t)| < \psi(t)$ for all $x \in A$ and $t>\rho(x)$.

The property of o-minimal structures that will play a crucial role in
this article concerns the two possibilities for asymptotic behavior
of functions $f: \BB R \to \BB R$ belonging to the structure.

\begin{df}
A structure on $(\BB R,+,\cdot)$ is {\em polynomially bounded} if for every
function $f: \BB R \to \BB R$ belonging to the structure, there exists some
$N \in \BB N$ (depending on $f$) such that $f(t) = O(t^N)$ as $t \to +\infty$.
A structure on $(\BB R,+,\cdot)$ is {\em exponential} if it contains $\exp$.
\end{df}

\begin{thm}[Growth dichotomy]
Either $\mathfrak S$ is polynomially bounded, or it is exponential.

If $\mathfrak S$ is polynomially bounded, then for every $f: \BB R \to \BB R$
belonging to $\mathfrak S$, either $f$ is ultimately identically equal to $0$,
or there exist $c,r \in \BB R$, $c \ne 0$, such that
$x \mapsto x^r : (0,\infty) \to \BB R$ belongs to $\mathfrak S$ and
$f(t) = c t^r + o(t^r)$ as $t \to +\infty$.
\end{thm}

\begin{cor}
The structures $\mathfrak S^{alg}$ and $\mathfrak S^{an}$ are polynomially
bounded.
\end{cor}

\begin{rem}
For concreteness, we state our results in the subanalytic setting.
That is, we assume that the symplectic manifold $M$, the group $G$ and
its action on $M$, the symplectic form $\omega$ and the moment map $\mu$
are subanalytic, i.e. belong to $\mathfrak S^{an}$.
This setting includes the vast majority of
examples of interest.
However, all results of this article hold in the setting of any o-minimal
structure which is polynomially bounded.
\end{rem}

Let $M$ be a closed ${\cal C}^1$ submanifold of $\BB R^d$ and
fix any norm $\|\,.\,\|_{\BB R^d}$ on $\BB R^d$.
%For $R \in \BB R$, let
%$$
%M_{\le R} = \{ x \in M ;\: \|x\|_{\BB R^d} \le R\}.
%$$
The manifold $M$ has a Riemannian metric induced by the standard metric on
$\BB R^d$.
%Hence
%$Vol \bigl( \{ x \in M ;\: \|x\|_{\BB R^d} \le R\} \bigr)$
%makes sense for $R \in \BB R$.
In light of the above results it is natural to expect

\begin{lem}  \label{polynomial1}
%Let $M$ be a closed ${\cal C}^1$ oriented submanifold of $\BB R^d$ and
%fix any norm $\|\,.\,\|_{\BB R^d}$ on $\BB R^d$.
%For $R \in \BB R$, let
%$$
%M_{\le R} = \{ x \in M ;\: \|x\|_{\BB R^d} \le R\}.
%$$
%The manifold $M$ has a Riemannian metric induced by the standard metric on
%$\BB R^d$. Hence $Vol(M_{\le R})$ makes sense.
If $M$ belongs to a polynomially bounded o-minimal structure $\mathfrak S$
(such as $\mathfrak S^{alg}$ or $\mathfrak S^{an}$), then there exists an
$N \in \BB N$ such that the function $f(R): \BB R \to \BB R$,
\begin{equation}  \label{funct1}
f(R)= Vol \bigl( \{ x \in M ;\: \|x\|_{\BB R^d} \le R\} \bigr),
\end{equation}
is $O(R^N)$, as $R \to +\infty$.

More generally, let $\alpha$ be a differential form on $M$ of top degree.
If both $M$ and $\alpha$ belong to a polynomially bounded o-minimal
structure $\mathfrak S$, then there exists an $N \in \BB N$ such that the
function $f(R): \BB R \to \BB R$,
\begin{equation}  \label{funct2}
f(R) = \int_{\{ x \in M ;\: \|x\|_{\BB R^d} \le R\}} |\alpha|,
\end{equation}
is $O(R^N)$, as $R \to +\infty$.
\end{lem}

\begin{lem}  \label{polynomial2}
Let $\mu: M \to \BB R$ be a proper function on a manifold $M$.
If $M$ and $\mu$ belong to a polynomially
bounded o-minimal structure $\mathfrak S$, then there exists an
$N \in \BB N$ such that the function $f(R): \BB R \to \BB R$,
\begin{equation}  \label{funct3}
f(R) = Vol \bigl( \{ x \in M ;\: \mu(x)=R \} \bigr),
\end{equation}
is $O(R^N)$, as $R \to +\infty$.

More generally, let $\alpha$ be a differential form on $M$.
If $M$, $\mu$ and $\alpha$ belong to a polynomially bounded o-minimal
structure $\mathfrak S$, then there exists an $N \in \BB N$ such that
the function $f(R): \BB R \to \BB R$,
\begin{equation}  \label{funct4}
f(R) = \int_{\{ x \in M ;\: \mu(x)=R \}} |\alpha|,
\end{equation}
is $O(R^N)$, as $R \to +\infty$.
\end{lem}

\begin{rem}
Note that the functions (\ref{funct1}), (\ref{funct2}), (\ref{funct3}) and
(\ref{funct4}) themselves need not belong to $\mathfrak S$.
Consider, for example, $M = \BB R$,
$$
\alpha(x) =
\begin{cases} dx & \text{if $|x| \le 1$;} \\
\frac{dx}{|x|} & \text{if $|x| \ge 1$.}
\end{cases}
$$
Both $M$ and $\alpha$ belong to $\mathfrak S^{min}$ and hence to any
polynomially bounded o-minimal structure $\mathfrak S$. However,
$$
f(R) = \int_{[-R,R]} \alpha =
\begin{cases}
0 & \text{if $R \le 0$;}  \\
2R & \text{if $0 < R \le 1$;}  \\
2 + 2\log R & \text {if $R>1$}
\end{cases}
$$
which cannot belong to $\mathfrak S$, since any structure containing
such a function $f(R)$ also contains the functions $\log$ and $\exp$.
\end{rem}

These two lemmas follow immediately from the Cell Decomposition Theorem
(Theorem \ref{cell}) described in the next subsection.
First we reduce the general case to the case when $M$ is a single cell,
then apply induction on the dimension of the cell.

\subsection
{Cells and cell decomposition}

In this subsection we continue to assume that $\mathfrak S$ is an o-minimal
structure on $(\BB R,+,\cdot)$. However, we make no assumptions on whether
$\mathfrak S$ is polynomially bounded or exponential.
Fix a positive integer $p$. We define the ${\cal C}^p$ cells in $\BB R^n$
as certain ${\cal C}^p$ submanifolds of $\BB R^n$ belonging to $\mathfrak S_n$.

\begin{df}  \label{cell-def}
Let $(i_1,\dots,i_n)$ be a sequence of zeroes and ones of length $n$.
An {\em $(i_1,\dots,i_n)$-cell of class ${\cal C}^p$} is a subset of $\BB R^n$
contained in $\mathfrak S$ obtained by induction on $n$ as follows:
\begin{enumerate}
\item
The ${\cal C}^p$ cells in $\BB R^1$ are just the singleton sets $\{r\}$
and the open intervals $(a,b)$, $-\infty \le a < b \le +\infty$.
The singletons $\{r\}$ are regarded as $(0)$-cells and the open intervals
$(a,b)$ as $(1)$-cells;
\item
Suppose $(i_1,\dots,i_n)$-cells of class ${\cal C}^p$ are already defined. 
Let $D \in \mathfrak S_n$ be a ${\cal C}^p$ $(i_1,\dots,i_n)$-cell, and let
$f: D \to \BB R$ of class ${\cal C}^p$ belong to $\mathfrak S$.
Then
$$
\operatorname{graph}(f) = \{ (x,r) \in D \times \BB R ;\: r=f(x) \}
$$
is an $(i_1,\dots,i_n,0)$-cell of class ${\cal C}^p$ in $\BB R^{n+1}$.
Let $g: D \to \BB R$ of class ${\cal C}^p$ be another function
contained in $\mathfrak S$ such that $f(x)<g(x)$ for all $x \in D$;
then the sets
\begin{align*}
& D \times \BB R  \\
& \{ (x,r) \in D \times \BB R ;\: r<f(x) \}  \\
& \{ (x,r) \in D \times \BB R ;\: r>f(x) \}  \\
& \{ (x,r) \in D \times \BB R ;\: f(x)< r < g(x) \}
\end{align*}
are $(i_1,\dots,i_n,1)$-cells of class ${\cal C}^p$ in $\BB R^{n+1}$.
\end{enumerate}
\end{df}

For example, a $(0,0)$-cell in $\BB R^2$ is a one point set,
a $(0,1)$-cell in $\BB R^2$ is a vertical interval,
a $(1,0)$-cell of class ${\cal C}^p$ in $\BB R^2$ is the graph of a
${\cal C}^p$ function defined on an interval and contained in $\mathfrak S$.
A $(1,\dots,1)$-cell in $\BB R^n$ is always open.
Each cell is connected and the dimension of an $(i_1,\dots,i_n)$-cell is
$i_1+\dots+i_n$.

Similarly, we define a {\em ${\cal C}^p$ cell decomposition} of $\BB R^n$ --
a special kind of partition of $\BB R^n$ into finitely many ${\cal C}^p$
cells.

\begin{df}
\begin{enumerate}
\item
A ${\cal C}^p$ cell decomposition of $\BB R$ is a collection of intervals and
points of the form
$$
\{ (-\infty, a_1), (a_1,a_2),\dots,(a_k,+\infty),\{a_1\},\dots,\{a_k\}\},
$$
with $a_1< \dots <a_k$ real numbers. (For $k=0$ this is just
$\{(-\infty,+\infty)\}$.)

\item
A ${\cal C}^p$ cell decomposition of $\BB R^{n+1}$ is a finite partition
${\cal D}$
of $\BB R^{n+1}$ into ${\cal C}^p$ cells such that the set of projections
$\{ \pi(D) ;\: D \in {\cal D} \}$ is a decomposition of $\BB R^n$, where
$\pi: \BB R^{n+1} \to \BB R^n$ is the projection on the first $n$ coordinates.
\end{enumerate}
\end{df}

Given a finite number of subsets $A_1,\dots,A_k \subset \BB R^n$,
we say that a partition ${\cal D}$ of $\BB R^n$ is
{\em compatible with $\{A_1,\dots,A_k\}$}
if for each $i$, $1\le i \le k$, and each $D \in {\cal D}$
either $D \subset A_i$ or $D \cap A_i = \varnothing$.

\begin{thm}[Cell decomposition]  \label{cell}
\begin{enumerate}
\item
Given $A_1,\dots,A_k \in \mathfrak S_n$, there is a ${\cal C}^p$ cell
decomposition of $\BB R^n$ compatible with $\{A_1,\dots,A_k\}$.
\item
For every function $f: A \to \BB R$ belonging to $\mathfrak S$,
$A \subset \BB R^n$ , there is a ${\cal C}^p$ decomposition ${\cal D}$ of
$\BB R^n$ compatible with $\{A\}$ such that the restriction
$f|_D: D \to \BB R$ is of class ${\cal C}^p$ for each
$D \in {\cal D}$ with $D \subset A$.
\end{enumerate}
\end{thm}

\noindent
{\bf Whitney stratification.}
Given $A_1,\dots,A_l \in \mathfrak S_n$, there is a finite ${\cal C}^p$
Whitney stratification of $\BB R^n$ compatible with $\{A_1,\dots,A_l\}$,
with each stratum a ${\cal C}^p$ cell in $\BB R^n$.

\separate

\separate

\noindent
{\em E-mail address:} {\tt mlibine@indiana.edu}

\noindent
{\em Department of Mathematics, Indiana University,
Rawles Hall, 831 East 3rd St, Bloomington, IN 47405}

\end{document}